\documentclass[graybox]{svmult}

\usepackage{mathptmx}       
\usepackage{helvet}         
\usepackage{courier}        
\usepackage{type1cm}        
\usepackage{makeidx}         
\usepackage{graphicx}        
\usepackage{multicol}        
\usepackage[bottom]{footmisc}

\newtheorem{thm}[subsection]{Theorem}
\newtheorem{prop}[subsection]{Proposition}
\newtheorem{lem}[subsection]{Lemma}
\newtheorem{cor}[subsection]{Corollary}

\usepackage[all]{xy}
\usepackage{amssymb}
\usepackage{amsmath}

\def\r{\succ}
\def\l{\prec}

\def\ss{\sigma}

\def\t{\otimes}

\newcommand{\KK}{\mathbb K}

\def\alg{\textrm{-}\mathsf{alg}}
\def\Sy{\mathbb{S}}

\newenvironment{proo}{\begin{trivlist} \item{\emph{Proof.}}}
  {\hfill $\square$ \end{trivlist}}

\def\ZZ{{\mathbb{Z}}}
\def\QQ{{\mathbb{Q}}}

\def\KK{{\mathbb{K}}}

\def\arbreA{\vcenter{\xymatrix@R=3pt@C=3pt{
&& \\
&*{}\ar@{-}[ur] \ar@{-}[ul] \ar@{-}[d]     &\\
&&
}}}

\def\arbreBA{\vcenter{\xymatrix@R=2pt@C=2pt{
&&&&\\
&&&*{}\ar@{-}[ul] & \\
&&*{}\ar@{-}[uurr] \ar@{-}[uull] \ar@{-}[d]     &&\\
&&&&
}}}

\def\arbreAB{\vcenter{\xymatrix@R=2pt@C=2pt{
&&&&\\
&*{}\ar@{-}[ur] &&& \\
&&*{}\ar@{-}[uurr] \ar@{-}[uull] \ar@{-}[d]     &&\\
&&&&
}}}

\def\arbreBB{\vcenter{\xymatrix@R=2pt@C=2pt{
&&*{}&&\\
&&&& \\
&&*{}\ar@{-}[uurr] \ar@{-}[uull] \ar@{-}[d] \ar@{-}[uu]     &&\\
&&&&
}}}

\def\arbreABC{\vcenter{\xymatrix@R=1pt@C=1pt{
&&&&&&\\
&*{}\ar@{-}[ur] &&&&& \\
&&*{}\ar@{-}[uurr] &&&&\\
&&&*{}\ar@{-}[uuurrr] \ar@{-}[uuulll] \ar@{-}[d] &&&\\
&&&&&&
}}}

\def\arbreBAC{\vcenter{\xymatrix@R=1pt@C=1pt{
&&&&&&\\
&&&*{}\ar@{-}[ul] &&& \\
&&*{}\ar@{-}[uurr] &&&&\\
&&&*{}\ar@{-}[uuurrr] \ar@{-}[uuulll] \ar@{-}[d] &&&\\
&&&&&&
}}}

\def\arbreACA{\vcenter{\xymatrix@R=1pt@C=1pt{
&&&&&&\\
&*{}\ar@{-}[ur] &&&&*{}\ar@{-}[ul] & \\
&&&&&&\\
&&&*{}\ar@{-}[uuurrr] \ar@{-}[uuulll] \ar@{-}[d] &&&\\
&&&&&&
}}}

\def\arbreCAB{\vcenter{\xymatrix@R=1pt@C=1pt{
&&&&&&\\
&&&*{}\ar@{-}[ur] &&& \\
&&&&*{}\ar@{-}[uull] &&\\
&&&*{}\ar@{-}[uuurrr] \ar@{-}[uuulll] \ar@{-}[d] &&&\\
&&&&&&
}}}

\def\arbreCBA{\vcenter{\xymatrix@R=1pt@C=1pt{
&&&&&&\\
&&&&&*{}\ar@{-}[ul] & \\
&&&&*{}\ar@{-}[uull] &&\\
&&&*{}\ar@{-}[uuurrr] \ar@{-}[uuulll] \ar@{-}[d] &&&\\
&&&&&&
}}}

\begin{document}

\title*{Exponential series without denominators}
\titlerunning{Exponential series without denominators} 
\author{Jean-Louis Loday}
\authorrunning{J.-L. Loday}
\institute{J.-L. Loday \at 
Institut de Recherche Math\'ematique Avanc\'ee,
    CNRS et Universit\'e de Strasbourg \\
    Zinbiel Institute of Mathematics\\ \email{loday@math.unistra.fr}}
\maketitle

\abstract{For a commutative algebra which comes from a Zinbiel algebra the exponential series can be written without denominators. When lifted to dendriform algebras this new series satisfies a functional equation analogous to the Baker-Campbell-Hausdorff formula. We make it explicit by showing that the obstruction series is the sum of the brace products. In the multilinear case we show that the role the Eulerian idempotent is played by the iterated pre-Lie product.}

\section*{Introduction} \label{S:int} The classical exponential series 
$$\exp(x)= 1+ x+ \frac{x^{2}}{2!}+\cdots + \frac{x^{n}}{n!}+ \cdots $$
can be written without denominators provided one assumes some properties on the commutative algebra. For instance in a Zinbiel algebra the  term $\frac{x^{n}}{n!}$ can be replaced by the iterated product $x^{\l n}:= x\l (x^{\l n-1})$, with $x^{\l 1}=x$. Recall that a Zinbiel algebra is a commutative algebra in which the product splits as $xy = x\l y + y\l x$. Matrices with coefficients in a commutative algebra (resp.\ Zinbiel algebra) are endowed with a structure of associative algebra (resp.\ dendriform algebra).  Recall that a dendriform algebra is an associative algebra whose product is the sum of two binary operations: $xy = x\l y + x\r y$,  supposed to  satisfy  some relations. There is a unique way to extend the exponential series from commutative algebras to associative algebras. However there are several ways to extend the exponential series from Zinbiel algebras to dendriform algebras: either take $\frac{x^{n}}{n!}$ or $x^{\l n}$ or $x^{\r n}$. These choices are equal in a Zinbiel algebra, but different in a dendriform algebra. 
We adopt the notation
$$e(x):=  1+ x+ x^{\l 2}+\cdots + x^{\l n}+ \cdots ,\quad e'(x):=  1+ x+ x^{\r 2}+\cdots + x^{\r n}+ \cdots $$
for the exponential series without denominators.
The classical functional equation 
$$\exp(x) \exp(y)=\exp(x+y)$$
 holds only when we are in the commutative or in the Zinbiel setting. In this last case it can be written:
 $$ e(x)e(y) = e(x+y).$$
 In the associative setting there is an error term for $\exp$ which is called the Baker-Campbell-Hausdorff formula. 
We prove that, in the dendriform context, the functional equation for the exponential without denominators is:

$$e(x)e(y) = e\big(x + e(x)\r y \l e'(-x)\big)\ ,$$
which can be written
$$e(x)e(y)= e(x+y+\sum_{n\geq 1} \{\underbrace{x, \ldots , x}_n;y\})\ ,$$
where $\{ -, \cdots ,- ; - \}$ is the brace product.
\medskip

In the associative framework the BCH-formula is best understood by looking at the multilinearized case. Then, the multilinear polynomial involved in the functional equation of the exponential turns out to be the Eulerian idempotent. In the dendriform case, we show that the analogue of the Eulerian idempotent for the exponential without denominators is played by the iterated pre-Lie product
$$h(x_{1},\ldots , x_{n})= \{x_{1}, \{ x_{2}, \{ \cdots \{ x_{n-1} , x_{n}\}\cdot\}\}\},$$
where the pre-Lie product is $\{x,y\} := x\r y - y\l x$.

\medskip

Here is the content of this paper.
In the first section we recall the notions of Zinbiel algebras, dendriform algebras and some of their properties. In the second section we introduce and study the exponential without denominators. 
In particular we compute its inverse,  the logarithm without denominators, which is closely related to the $1\frac{1}{2}$-logarithm of Kontsevich.
In section 3 we prove the BCH-type formula for the exponential series without denominators. As a corollary we show that the multilinear obstruction to the additivity of the exponential is the iterated pre-Lie product.
\medskip

In this paper $\KK$ denotes a unital commutative ring (for instance $\ZZ$) which is the ground ring. We sometimes need to suppose that $\KK$ contains $\QQ$. The tensor product over $\KK$ is simply denoted by $\t$. 
\bigskip


\section{Zinbiel and dendriform algebras}\label{Zinbielalg} We introduce the notion of Zinbiel algebra (called dual Leibniz algebra in \cite{Loday95b}) and its relationship with commutative algebras. We also introduce the notion of dendriform algebras because matrices over a Zinbiel algebra bear the structure of a dendriform algebra.

\subsection{Definition}\label{Zinbielalg} A \emph{Zinbiel algebra} is a module $A$ over $\KK$ equipped with a binary operation $x\l y$ such that
$$(x\l y)\l z = x\l (y\l z+z\l y).$$
Symmetrizing the Zinbiel operation, that is defining $xy:=x\l y+y\l x$, we get a product on $A$ which is commutative and associative. Indeed, one gets
\begin{multline*}
(xy)z= (x\l y +y\l x)\l z + z\l (xy)= x\l (yz)+ y\l (xz) + z\l (xy)=\\ x\l (yz)+ y\l (zx) + z\l (yx)=x\l (yz)+ (y\l z)\l x+(z\l y)\l x = x(yz).
\end{multline*}
Hence there is a forgetful functor between categories of algebras:
$$Zinb\alg \to Com\alg.$$

\subsection{Free Zinbiel algebra} It is shown in \cite{Loday01} that the free Zinbiel algebra over the vector space $V$ is the reduced tensor module
$$\bar T(V):= V\oplus \cdots \oplus V^{\t n} \oplus \cdots $$
where the generic element $v_1\cdots v_n\in V^{\t n}$ corresponds to the product $v_1\l (v_2\l ( \cdots \l v_n))$. As a consequence the space of $n$-ary operations of the operad $Zinb$ is the $\Sy_n$-module $Zinb(n)=\KK[\Sy_n]$. More precisely we have:

\begin{thm} Let $V$ be a vector space and let $\bar T(V)$ be the reduced tensor module over $V$. The half-shuffle:
$$v_{1}\cdots v_{p}\l v_{p+1}\cdots v_{p+q}:= \sum_{\ss} v_{1} v_{\ss^{-1}(2)}\cdots v_{\ss^{-1}(p+q)}$$
where $\ss$ is a $(p-1,q)$-shuffle acting on $\{2, \ldots, p+q\}$, makes $\bar T(V)$ into the free Zinbiel algebra on $V$. The associated commutative algebra is the (nonunital) shuffle algebra, denoted by $\bar T^{sh}(V)$. The Zinbiel algebra $(\bar T (V), \l)$ is in fact the free Zinbiel algebra over $V$, denoted by $Zinb(V)$.
\end{thm}
\begin{proo} See \cite{Loday01}.
\end{proo}

\subsection{Example} If $V$ is one-dimensional spanned by $x$, then $Zinb(\KK\, x)$ is spanned by $x^{\l n}$ for $n\geq 1$. The operations are given by
$$x^{\l p} \l x^{\l q} = \binom{p+q-1}{q-1} x^{\l p+q} \quad \mathrm{and} \quad x^{\l p}x^{\l q}= \binom{p+q}{q} x^{\l p+q}\ .$$

\subsection{Dendriform algebra} A \emph{dendriform algebra} is a $\KK$-module $A$ equipped with two linear maps (binary operations)
$$\l : A\otimes A \to A \quad \mathrm{and } \quad \r : A\otimes A \to A $$
called the \emph{left operation} and the \emph{right operation} respectively, satisfying the following three relations
\begin{displaymath}
\left\{\begin{array}{rcl}
(x\l y)\l z &=& x\l (y\l z)+ x\l (y\r z)   , \\
(x\r y)\l z &=& x\r (y\l z) , \\
(x\l y)\r z + (x\r y)\r z&=& x\r (y\r z).
\end{array}
\right .
\end{displaymath}

From these axioms it follows readily that the binary operation
$$xy := x\l y + x \r y\ .$$
is associative. Under this notation the first relation becomes 
$$(x\l y)\l z = x\l (yz).$$
 In the proof of the main Theorem we allow ourselves to write  $x\l yz$ in place of $x\l (yz)$.

Let us mention that numerous combinatorial Hopf algebras come with a dendri\-form structure, cf.\ for instance \cite{Aguiar04, BurgunderRonco10, LodayRonco10, NovelliThibon07}.
 
\subsection{Commutative dendriform algebra} By definition a commutative dendriform algebra is a dendriform algebra whose left and right operations are related by the following symmetry relation
$$ x\l y = y \r x, \textrm{ for any } x \textrm{ and } y .$$

By direct inspection we see that a Zinbiel algebra is a commutative dendriform algebra and vice-versa (cf.\ \cite{Loday01}).

\begin{prop} The module of $n\times n$-matrices $\mathcal{M}_{n}(A)$ with coefficients in the Zinbiel algebra $A$ is a dendriform algebra.
\end{prop}

\begin{proo} It is straightforward to check that the matrices over a dendriform algebra is still a dendriform algebra. The formulas are like in the classical case. Since a Zinbiel algebra is a particular case of a dendriform algebra, we are done.
\end{proo}

\subsection{Dendriform calculus} We recall some results from \cite{Loday01} about  computation in a dendriform algebra $A$.

It is often helpful to suppose that $A$ is equipped with a unit $1$. So $A=\KK\, 1 \oplus I$ where $I$ is a dendriform algebra, and the relationship between $1$ and the elements of $I$ are given by:
\begin{align*}
x\l 1 = x, &\qquad  1\l x = 0, \\
x\r 1 = 0, &\qquad  1\r x = x. 
\end{align*}
Observe that $1\l 1$ and $1\r 1$ are not defined, so, in order for an expression like $(1+u)(1+v)$ to make sense, we have to write it as 
$$  (1+u)(1+v)= 1+u+v+uv=1+u+v+u\l v+ u\r v \ .$$

The \emph{free unital dendriform algebra} on one generator is spanned by the planar rooted binary trees. We denote by $PBT_{n}$ the set of planar rooted binary trees $t$ with $n$ leaves. In low dimension we have:
\begin{displaymath}
PBT_1 = \{\  \vert \ \}\  , \quad PBT_2=\big\{ \arbreA\big\}\  , \quad 
PBT_3=\Big\{\arbreBA\ ,\  \arbreAB \Big\}\ ,
\end{displaymath}
\begin{displaymath}
{PBT_4=\bigg\{ \arbreABC}\ ,\  {{}\arbreBAC },\  {{}\arbreACA },\  {{}\arbreCAB },\  {{}\arbreCBA }\bigg\}\ .
\end{displaymath}

When $t$ has $n+1$ leaves (i.e.\ $t\in PBT_{n+1}$), it determines an $n$-ary operation. Applied to the generic $n$-tuple $x_{1} \cdots x_{n}$ we denote the result by $t(x_{1} \cdots x_{n})$. For instance:
$$|()= 1, \quad \arbreA(x) = x, \arbreAB(xy)= x\r y, \quad \arbreBA(xy)= x\l y.$$
More generally, if the tree $t$ is the grafting of $t^{l}$ and $t^{r}$, denoted by $t= t^{l}\vee t^{r}$, then 
$$t(x_{1} \cdots x_{n}) = t^{l}(x_{1} \cdots x_{i-1})\r x_{i}\l  t^{r}(x_{i+1} \cdots x_{n}).$$

The free dendriform algebra over the vector space $V$ is $Dend(V)= \bigoplus_{n} Dend_{n}\t V^{\t n}$  where $Dend_{n}= \KK[PBT_{n+1}]$. In order to check an equality of multilinear elements in the free dendriform algebra over the set $\{x_{1}, \ldots , x_{n}\}$ it is equivalent to check the equality for each permutation of the variables individually (in other terms the operad $Dend$ is a nonsymmetric operad).

It is sometimes helpful to adopt the notation  $x^{t}:=  t(x\cdots x)$ or simply $t$ when there is only one variable into play. For instance we have:
$x^{t}= x^{\l n}= x\l (x\l ( \cdots \l x))$ for $t$ the right comb with $n+1$ leaves. It can be shown that 
$$x^{n}= \sum _{t\in PBT_{n+1}} x^{t}$$
in any dendriform algebra.

We already introduced the notation $x^{\l n}$. Similarly we define $x^{\r n}:=(x^{\r n-1})\r x $, with $x^{\r 1}:=x$.

If we now work in a Zinbiel algebra, then it can be shown that $x^{t}= \# \varphi^{-1}(t)\, x^{\l n}$ where $\varphi : \Sy_{n}\to PBT_{n+1}$ is the surjective map constructed in loc.cit. As a consequence we have
$$x^{n}= \sum _{t\in PBT_{n+1}} x^{t}=  \sum _{t\in PBT_{n+1}}  \# \varphi^{-1}(t)\,  x^{\l n}=  \sum _{\ss\in \Sy_{n}} x^{\l n}= n!\, x^{\l n},$$
as we already know.

\subsection{Pre-Lie product and brace products} Recall that a binary operation is said to be \emph{left pre-Lie} if its associator is symmetric in the first two variables. In a dendriform algebra the binary operation 
$$\{x,y\}:= x\r y - y\l x$$
 is a left pre-Lie product. It is a direct consequence of the dendriform axioms. More generally, one can form the following $(n+1)$-ary operation:
$$\{x_1, \ldots, x_n ; y\} := \sum_{i=0}^{n} (-1)^{n-i}(x_1 \l ( x_2 \l ( \cdots \l x_i ))) \r y \l (((x_{i+1}\r \cdots )\r x_{n-1})\r x_n
)\ .$$
In \cite{Ronco00} M.\ Ronco showed that these operations are primitive for the Hopf structure of the free dendriform algebra and that they satisfy the axioms of a \emph{brace algebra}. They are called \emph{brace products}. Observe that $\{x;y\}=\{x,y\}$ is the left pre-Lie product.

\subsection{A commutative diagram of algebras}
The relationship between all these types of algebras can be summarized by the existence of a commutative diagram of categories of algebras  \cite{Loday01, Ronco00}:

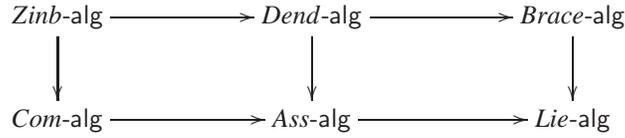
$$\xymatrix{
 Zinb\alg\ar[rr] \ar[d] & &Dend\alg \ar[d] \ar[rr] &&Brace\alg\ar[d] \\
Com\alg\ar[rr]  & &Ass\alg\ar[rr] &&Lie\alg
}$$

Considering free algebras on one generator $x$ we get the polynomials $\KK[x]$ for $Zinb, Com$ and $Ass$, and we get $\bigoplus_{n} \KK[PBT_{n+1}]$ for $Dend$:

$$
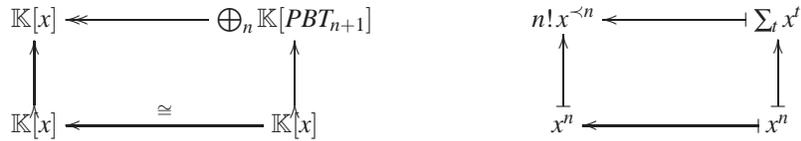
\xymatrix{
 \KK[x]& & \bigoplus_{n}  \KK[PBT_{n+1}]\ar@{->>}[ll] &&n!\, x^{\l n} & &\sum_{t} x^{t}\ar@{|->}[ll] \\
\KK[x]  \ar@{>->}[u] & & \KK[x]      \ar[ll]_{\cong} \ar@{>->}[u] &&x^{n}\ar@{|->}[u]& & x^{n} \ar@{|->}[u]\ar@{|->}[ll] 
}$$

Observe that there are various possibilities to lift the element $\frac{x^n}{n!}$ to the free dendrifrom algebra, for instance $\frac{x^t}{\# \varphi^{-1}(t)}$ for any $t\in PBT_{n+1}$.


\section{Exponential series}\label{expseries}
We introduce the exponential series without denominators and we compute their inverse for the associative product and for composition (logarithm without denominators).

\subsection{Exponential series without denominators} Let $V$ be a module over $\KK$ and let $Dend(V)^{\wedge}$ be the infinite product
$$Dend(V)^{\wedge}:= \prod_{n} Dend_{n}\t V^{\t n}.$$
For any $x\in V$ a series in $x$ is an element of $Dend(V)^{\wedge}$ made out of products of the only variable $x\in V$. 

  There are several different ways to extend the exponential series $\exp(x) = e(x)= e'(x)$ on Zinbiel algebras to dendriform algebras. For instance we have
  
\begin{enumerate}

\item $\exp(x)= 1+ x+ \frac{x^{2}}{2!}+\cdots + \frac{x^{n}}{n!}+ \cdots $, the classical exponential series when $\QQ\subset \KK$,

\item $e(x):= 1+ \sum_{n\geq 1}x^{\l n},\  e'(x):=1+ \sum_{n\geq 1}x^{\r n}$, the exponential series without denominators,

\item $ee(x):= 1+ \sum_{n\geq 1}\frac{1}{2}( x^{\l n}+ x^{\r n})$, when $2$ is invertible in $\KK$.
\end{enumerate}

Observe that we have $e(x) = 1 + x\l e(x)$.

\begin{lem}\label{lemmecle} In the dendriform context we have:

$$e'(-x)e(x)=1= e(x)e'(-x).$$
\end{lem}
\begin{proo} Let us first prove that $x^{\r i}\r x^{\l j}= x^{\r i+1}\l x^{\l  j-1}$. We have
\begin{align*}
x^{\r i}\r x^{\l j}&= x^{\r i}\r (x\l x^{\l j-1})\\
&= (x^{\r i}\r x)\l x^{\l j-1}\\
&= x^{\r i+1}\l x^{\l j-1}\ .\\
\end{align*}
Expanding the product $ e'(-x)e(x)$ we get 
$$x^{\l n} + \cdots + (-1)^{j}( x^{\r i}\r x^{\l j} + x^{\r i}\l x^{\l j}) +\cdots + (-1)^{n}x^{\l n}$$
in degree $n>0$. This element is $0$ as a consequence of the preceding formula. The second formula is an immediate consequence of the first one.
\end{proo}

 \begin{prop} In the dendriform context, the inverse of the exponential series $E(x):=e(x)-1$ for composition is the series
$$L(x):= x\l (1-x +x^{2}- x^{3} + \cdots +(-1)^{n} x^{n}+\cdots) \quad ,$$
called the \emph{logarithm without denominators}.
\end{prop}

\begin{proo} Let us write $L(x)= x + \varphi_{2}(x)+   \cdots +\varphi_{n}(x)+\cdots$ for the inverse of $E(x)$. Since $E(x) = x\l (1+ E(x))$, by replacing $x$ by $L(x)$ in this equality we get
$$x = L(x) \l (1+ x).$$
Hence we have $\varphi_{2}(x)= -x\l x$ and $\varphi_{n+1}(x)+ \varphi_{n}(x)\l x= 0$. By induction we suppose that $ \varphi_{n}(x)= (-1)^{n-1}x\l x^{n-1}$. We compute:
\begin{align*}
\varphi_{n+1}(x)&= -\varphi_{n}(x)\l x,\\
 &= - (-1)^{n-1}(x\l x^{n-1})\l x,\\
&= (-1)^{n}x\l (x^{n-1}\l x + x \r x^{n-1}),\\
&=  (-1)^{n}x\l (x^{n}).
\end{align*}
\end{proo}

\subsection{Remark} This proposition gives a quick proof of the fact that $\exp$ and $\log$ are inverse to each other. This proof is similar to the proof which uses the integral definition of the logarithm: $\log(1+x)= \int \frac{dx}{1+x}$.

\subsection{Zinbiel algebras in characteristic $p$} Let us suppose that $\KK$ is a characteristic $p$ field and let us work in the Zinbiel framework. Since $x^{n}= n!\, x^{{\l n}}$ the logarithm becomes a polynomial which can be written
$$\log(1+x) = L(x) := \sum_{i=1}^{p-1}(-1)^{i-1}\frac{x^{i}}{i} + (-1)^{p-1}(p-1)!\,  x^{\l p}.$$
First, the element $x^{\l p}$ is a divided power, that is $(A, \gamma(x):= x^{\l p})$ is a divided power algebra, cf.\ for instance \cite{Dokas10}. Second, the first part of the logarithm is the so-called ``one and a half logarithm'' introduced by Maxim Kontsevich in \cite{Kontsevich02}. In other words we have the following result.

\begin{prop} In a characteristic $p$ divided power algebra (resp.\ Zinbiel algebra) the one and a half logarithm plus $(-1)^{p-1}(p-1)!\, \gamma(x)$ is an invertible series whose inverse is the exponential series. \hfill $\square$
\end{prop}

\section{BCH-type formula for the exponential series without denominators}

In a Zinbiel algebra we have the functional equation $ e(x)e(y) = e(x+y)$. But, in a dendriform algebra, $x\r y \neq y\l x$ (even for $x=y$), hence there are nontrivial dendriform polynomials $H_{n}(x,y)$ in $x$ and $y$ such that
$$e(x)e(y) = e(x+y+ \cdots + H_{n}(x,y)+\cdots ).$$
For instance we obviously have $H_{2}(x,y)= x\r y - y\l x$, which is the (left) pre-Lie product $\{x,y\}$. Our aim is to compute $H_{n}(x,y)$, in fact to show that it is a brace product. 

\begin{thm}\label{mainthm} In a dendriform algebra the following equalities hold:

$$e(x)e(y) = e\big(x + e(x)\r y \l e'(-x)\big)\ ,$$

$$e'(x)e'(y) = e'\big( e'(-y)\r x \l e'(y)+y\big)\ ,$$

for  $e(x):= 1+ \sum_{n\geq 1}x^{\l n}$ and $ e'(x):=1+ \sum_{n\geq 1}x^{\r n}$.
\end{thm}

\begin{cor}\label{BCHnew}  The dendriform polynomial $H_{n}$ is the brace product:
$$H_{n}(x,y)= \sum_{i+j=n-1}(-1)^{j} (x^{\l i})\r y \l (x^{\r j})= \{\underbrace{x, \ldots, x}_{n-1} ; y\}\ .$$
\end{cor}

\begin{proo} [of Theorem \ref{mainthm}] We prove the first relation. Let us define
$$R(x,y):= e(x)e(y) - e\big(x + e(x)\r y \l e'(-x)\big)\ .$$
The aim is to prove that $R(x,y)=0$.

We use the following two relations: $e(z) = 1 + z\l e(z)$ and $e'(-x) e(x) =1$ (Lemma \ref{lemmecle}). On one hand we have
\begin{eqnarray*}
e\big(x + e(x)\r y \l e'(-x)\big)&=&1+ \big(x+e(x)\r y \l e'(-x)\big)\l (e(x)e(y)-R(x,y))\\
&=& 1+x\l(e(x)e(y))+ (e(x)\r y)\l (e'(-x)e(x)e(y))\\
&&- \Phi(x,y),
\end{eqnarray*}
where $\Phi(x,y):= (x + e(x)\r y \l e'(-x))\l R(x,y)$.
Then we have
\begin{eqnarray*}
e\big(x + e(x)\r y \l e'(-x)\big)&=& 1+x\l (e(x)e(y)) + e(x)\r y\l e(y)- \Phi(x,y)\\
&=& 1+x\l (e(x)e(y)) + e(x)\r E(y)- \Phi(x,y),
\end{eqnarray*}
by Lemma \ref{lemmecle}. 
On the other hand we have:
\begin{eqnarray*}
e(x)e(y) &=& e(x)(1+E(y)) \\
&=& e(x)+ e(x)\l E(y) +e(x)\r E(y). \\
\end{eqnarray*}
Hence we compute 
\begin{eqnarray*}
R(x,y)&=& e(x) +e(x)\l E(y))-1-x\l (e(x)e(y))+ \Phi(x,y)\\
&=& e(x) +e(x)\l E(y)-1-x\l (e(x)(1+E(y))+ \Phi(x,y)\\
&=&e(x) +e(x)\l E(y) -1-x\l e(x) - x\l (e(x)E(y)) + \Phi(x,y)\\
&=&e(x)\l E(y)  -x\l (e(x)E(y))+ \Phi(x,y)\\
&=& e(x)\l E(y) -(x\l e(x))\l E(y)+  \Phi(x,y)\\
&=& 1\l E(y) + \Phi(x,y)\\
&=& \Phi(x,y)\ .
\end{eqnarray*}
So we have proved that $R(x,y)$ satisfies the functional equation:
$$R(x,y)= (x+e(x)\r y \l e'(-x))\l R(x,y)\ .$$
Since $R(0,0)=0$ we see from this equation that the degree 1 part of $R(x,y)$ is also 0, and, by induction, the degree $n$ part of $R(x,y)$ is $0$ for any $n$. So $R(x,y)=0$ and we are done.

The second formula follows from the fact that the involution $\tau :Dend(V)\to Dend(V)$ which sends $\l$ to $\r$, $\r$ to $\l$ and $v_1 \cdots v_n$ to $v_n\cdots v_1$ is an isomorphism of dendriform algebras. 

The  Corollary is an immediate consequence.
\end{proo}

\subsection{Functional equation in one variable}  In the associative case, if $x=y$, then the functional equation of the classical exponential series is the same as in the commutative case. However in the dendriform case,  if $x=y$, then the functional equation  is different from the Zinbiel case because $x\l x\neq x\r x$. As an immediate corollary of the main result the term $H_{n}(x,x)$ in the formula 
$$e(x)e(x) = e\Big(2x + \sum_{n\geq 2}H_{n}(x,x)\Big)$$
is as follows. Let us denote by $lc^{i}$ (resp.\ $rc^{i}$) the left comb (resp.\ right comb) with $i$ leaves.

\begin{cor} In the dendriform context we have:
$$H_{n}(x,x) = \sum_{i+j = n+1\atop i\geq 1, j\geq 1} (-1)^{j} rc^{i} \vee lc^{j}\in \KK[PBT_{n+1}] \  for\  n\geq 2.$$
\end{cor}
\begin{proo} It follows from Corollary \ref{BCHnew} and from the following equalities in the free  dendriform algebra on one generator:
$$
rc^{i} (x, \ldots, x)= x^{\l i-1} ,\quad
lc^{i} (x, \ldots, x)= x^{\r i-1} .
$$
\end{proo}

\subsection{The multilinear case} For the classical exponential series in the associative case, the Baker-Campbell-Hausdorff formula takes the form 
$$\exp(x)\exp(y) = \exp(x+y + \cdots +BCH_m(x,y)+ \cdots )\ .$$
One of the way to compute the Lie polynomial $BCH_m(x,y)$ is to consider the multilinear version
$$\exp(x_1)\cdots \exp(x_n) = \exp(x_1+\cdots +x_n  + \cdots +BCH_m(x_1,\ldots , x_n)+ \cdots )\ .$$
It is a result of Dynkin that $BCH_m(x,y)$ can be computed out of the multilinear part of $BCH_n(x_1,\ldots , x_n)$ denoted $eul(x_1,\ldots , x_n)\in \QQ[\Sy_n]$. It is known, cf.\ for instance \cite{Loday94}, that $eul(x_1,\ldots , x_n)$ is the Eulerian idempotent.

Let us multilinearize similarly  the functional equation of the series $e(x)$ in the dendriform case:
$$e(x_1)\cdots e(x_n) = e(x_1+\cdots +x_n  + \cdots +H_m(x_1,\ldots , x_n)+ \cdots )\ .$$
 and denote by  $h(x_1,\ldots , x_n)$ the multilinear part of $H_n(x_1,\ldots , x_n)$.
 
 \begin{prop} The dendriform polynomial  $h(x_1,\ldots , x_n)$ is the iterated pre-Lie product:
$$h(x_{1},\ldots , x_{n})= \{x_{1}, \{ x_{2}, \{ \cdots \{ x_{n-1} , x_{n}\}\cdot\}\}\},$$
where $\{x,y\} := x\r y - y\l x$.
 \end{prop} 
 
 \begin{proo} In order to compute $e(x_1)\cdots e(x_n) $ we can apply the functional equation of Theorem \ref{mainthm} iteratively. It comes immediately that the polynomial 
 $h(x_1,\ldots , x_n)$ is the multilinear part of
 $$e(x_1)\r\Big(e(x_2)\r\big( \cdots (e(x_{n-1})\r x_n \l e'(-x_{n-1})) \cdots \big)\l e'(-x_{2})\Big)\l e'(-x_{1}).$$
 It comes
 $$h(x_1,\ldots , x_n)= \sum (-1)^{\ell}x_{i_1}\r \Big(x_{i_2} \r \big( \cdots (x_{i_k}\r x_n\l x_{j_\ell})\l \cdots  \big)\l x_{j_2}\Big)\l x_{j_1} ,$$
 where $\{ i_{1}\ \cdots\ i_{k}\ |\ j_1\ \cdots\ j_{\ell}\}$ is a $(k,l)$-shuffle of $\{1, \ldots , n-1\}$.
 
 It is clear that for $n=2$ we get $h(x_1, x_2)=x_1\r x_2 - x_2 \l x_1 = \{x_1,x_2\}$. Then, by induction we get the expected result.
  \end{proo}

\end{document}